\documentclass[reqno,12pt]{amsart}%,draft
\usepackage[mathscr]{eucal}
\usepackage{amsmath}
\usepackage{amssymb}
\usepackage{latexsym}
\usepackage{amsthm}
\usepackage{color}
\newtheorem{thm}{\indent Theorem}[section]

\newtheorem{lem}[thm]{\indent Lemma}
\newtheorem{prop}[thm]{\indent Proposition}

\newtheorem{rmk}{{\indent\bf Remark}}[section]

\newtheorem*{conjecture}{\indent Conjecture}

\newcommand{\lmx}{\left(\begin{matrix}}
\newcommand{\rmx}{\end{matrix}\right)}
\newcommand{\ldt}{\left|\begin{matrix}}
\newcommand{\rdt}{\end{matrix}\right|}

\newcommand{\lagl}{\langle}
\newcommand{\ragl}{\rangle}

\newcommand{\bbr}{{\mathbb R}}
\newcommand{\bbh}{{\mathbb H}}

\newcommand{\bbc}{{\mathbb C}}
\newcommand{\be}{\begin{equation}}
\newcommand{\ee}{\end{equation}}

\newcommand{\vfi}{\varphi}

\makeatletter
\numberwithin{equation}{section}
\begin{document}

\title [A sharp Bernstein theorem for affine maximal type surfaces]{{A sharp Bernstein theorem for affine maximal type surfaces}}

\author{Yalin Sun}
\address{School of Mathematics and Statistics,
Henan Normal University, Xinxiang 453007,
People's Republic of China}
\email{ylsun2024@126.com}

\author{Cheng Xing}
\address{
School of Mathematics and Statistics,
Henan Normal University, Xinxiang 453007,
People's Republic of China}
\email{xingcheng@htu.edu.cn}

\author{Ruiwei Xu}
\address{
School of Mathematics and Statistics,
Henan Normal University, Xinxiang 453007,
People's Republic of China}
\email{rwxu@htu.edu.cn}


\thanks{R. Xu is the corresponding author}

\date{}


\begin{abstract}
We investigate Bernstein theorem for the affine maximal type equation
\begin{equation*}
\sum_{i,j=1}^nf^{ij}(D^a)_{ij}=0,\qquad
D=\det(f_{ij}),\qquad a\neq0,
\end{equation*}
under completeness of the Calabi metric.
In dimension two, we prove that every smooth strictly convex solution
on a convex domain is a quadratic polynomial whenever
$a\notin[-\tfrac23,-\tfrac13]$, resolving a conjecture of Sun and Xu \cite{SX}.
For every $n\geq2$, we also construct explicit non-quadratic
Calabi complete solutions on convex domains for all
$a\in[-\tfrac n{n+1},-\tfrac1{n+1}]$.
These examples show that the rigidity range is sharp in dimension two.
\end{abstract}

\maketitle
\renewcommand{\baselinestretch}{1}
\renewcommand{\arraystretch}{1.2}
\catcode`@=11 \@addtoreset{equation}{section} \catcode`@=12
\makeatletter
\renewcommand{\theequation}{\thesection.\arabic{equation}}
\@addtoreset{equation}{section} \makeatother

%%%%%%%%%%%%%%%%
\section{Introduction}
The rigidity of complete parabolic affine hyperspheres is an important problem in affine differential geometry,  see \cite{Ca, CY, G, J, LSZH, LXSJ, P}. 
The main goal of this paper is the Bernstein problem for affine maximal type hypersurfaces, a natural generalization of parabolic affine hyperspheres.
Let $f$ be a smooth strictly convex solutions on a convex domain
$\Omega\subset\bbr^n$, $n\geq2$, of the fourth-order equation with a real constant exponent $a$:
\be\label{1.1}
\sum_{i,j=1}^nf^{ij}(D^a)_{ij}=0,\qquad D:=\det(f_{ij}),\qquad a\neq0,
\ee
where $f_{ij}=\tfrac{\partial^2f}{\partial x_i\partial x_j}$ and $(f^{ij})=(f_{ij})^{-1}$.
Throughout this paper, strict convexity means that the Hessian
matrix $(D^2f)=(f_{ij})$ is positive definite.
Following \cite{D}, we call (\ref{1.1}) the affine maximal type equation.
It is the Euler-Lagrange equation of the affine area functional \cite{D2,LXSJ}:
\[
\mathcal A(f,\Omega):=\int_\Omega D^{a+1}\,dx
\qquad (a\neq-1),
\]
with the logarithmic variant
\[
\mathcal A(f,\Omega):=\int_\Omega\ln D\,dx
\]
for the exceptional value \(a=-1\).

Several important equations arise as special cases of (\ref{1.1}).
For $a=-1$, it becomes the \emph{Abreu equation}, which appears in the differential geometry of toric varieties \cite{A, D-2}.
For $a=-\tfrac12$, it characterizes Calabi affine
extremal hypersurfaces \cite{XL22}.
In dimension two, these are also called \emph{Calabi affine maximal surfaces} \cite{SXX26}.
For the special case $a=-\tfrac{n+1}{n+2}$,  the PDE
(\ref{1.1}) characterizes \emph{affine maximal hypersurfaces} \cite{Ca,LXSJ}.
The graph
\[
M_f:=\{(x,f(x)):x\in\Omega\}
\]
has Blaschke metric
\[
G^B_{ij}:=\rho\, f_{ij},\qquad
\rho:=D^{-\tfrac{1}{n+2}}.
\]
Equation (\ref{1.1}) is equivalent to
\begin{equation}\label{1.2}
\Delta_M \rho=0,
\end{equation}
where \(\Delta_M\) is the Laplace-Beltrami operator of $G^B$. Thus \(M_\Omega\) is affine maximal if and only if \(\rho\) is harmonic with respect to the Blaschke metric.

In dimension two, Donaldson \cite{Do} developed the local theory of (\ref{1.1}) by introducing a change of variables that, near
points where $dD\neq0$, relates its solutions to pairs of solutions of a second-order linear PDE.
In dimension two, Chern's conjecture asserts that every entire
smooth strictly convex affine maximal graph is an elliptic
paraboloid \cite{Ch}.
Calabi's conjecture asserts the same conclusion under completeness
of the Blaschke metric \cite{Ca,LXSJ}.
In 2000, Trudinger and Wang \cite{TW} proved Chern's conjecture and proposed a higher-dimensional generalization under Euclidean
completeness.
For equation (\ref{1.1}), they  \cite{TW1} also proved that
entire locally uniformly convex solutions on $\bbr^2$ are quadratic when $a>0$.
Li and Jia \cite{LJ} established the same conclusion for $a\leq-\tfrac34$, providing a new analytic
proof of Chern's conjecture; see also \cite{Z} for a different proof when $a\in[-1,-\tfrac34)$.
For $n\geq2$, Jia and Li \cite{JL} obtained Bernstein results for entire solutions when
$|a|$ is sufficiently large; see also \cite{XLL}.
Further results can be found in \cite{CSL,D2,L-J-4,LXSJ,LXSJ1,XXS}.

On the other hand, non-quadratic complete solutions have also been constructed in recent years.
Du \cite{D} obtained Euclidean complete examples for $n\geq2$ and
$a\in[-\tfrac{n-1}{n},0)$.
Sun and Xu \cite{SX} constructed examples at
$a=-\tfrac n{n+1}$ and $a=-\tfrac1{n+1}$ that are complete
with respect to both the Euclidean and Calabi metrics. Recently,  Sun, Xing and Xu \cite{SXXE26} extended the parameter range of
Euclidean complete examples to $[-\tfrac n{n+1},0)$.

Another natural completeness assumption concerns the Calabi metric
\be\label{1.3}
G:=\sum_{i,j=1}^n f_{ij}dx_idx_j,
\ee
which Calabi \cite{Ca} originally called the {\em Schwarz-Pick metric}.
Li and Jia \cite{L-J-4} proved that every Calabi complete affine maximal hypersurface is an elliptic
paraboloid when $n=2$ or $n=3$. For the affine maximal type equation, the following more general
result holds; see Theorem 5.8.1 in \cite{LXSJ} or (2)  of  Theorem 1.1 in \cite{D2}.

\begin{thm}[\cite{LXSJ}] \label{1.0.0}
Let $f$ be a smooth strictly convex solution of the PDE (\ref{1.1}) on a convex domain  $\Omega\subset\bbr^n$, $n\geq2$. If the
Calabi metric (\ref{1.3}) is complete and
\be\label{1.4}
a\notin\left[-\tfrac{1}{2}-\tfrac{n-1}{4\sqrt{n}},-\tfrac{1}{2}+\tfrac{n-1}{4\sqrt{n}}\right],
\ee
then $f$ is a quadratic polynomial.
\end{thm}

Non-quadratic Calabi complete solutions are known for certain
exponents outside the range covered by Theorem 1.1.
Xu and Li \cite{XL22} first constructed such solutions
at $a=-\tfrac12$ in every dimension $n\geq2$ (see also \cite{XL21}).
Sun and Xu \cite{SX} investigated a special class of Calabi
hypersurfaces with negative constant sectional curvature and obtained
Calabi complete examples at $a=-\tfrac{n}{n+1}$ and $a=-\tfrac1{n+1}$.
They observed that in dimension two these exponents reduce to
$-\tfrac23$ and $-\tfrac13$, which differ from the corresponding
endpoints $-\tfrac12\pm\tfrac1{4\sqrt2}$ of Theorem 1.1 by only
$\tfrac1{4\sqrt2}-\tfrac16\approx0.0101$.
In view of this narrow gap, they conjectured that the rigidity
conclusion of Theorem 1.1 holds for all
$a\notin[-\tfrac23,-\tfrac13]$.

\begin{conjecture}[\cite{SX}]
Let $f$ be a smooth strictly convex function defined on a convex domain $\Omega$ in $\bbr^2$. If $f$ satisfies
the PDE (\ref{1.1}) with $a\notin[-\tfrac{2}{3},-\tfrac{1}{3}]$ and the Calabi metric is complete, then $f$ is
a quadratic polynomial.
\end{conjecture}

We prove this conjecture using the normalized Hessian map
\[
\mathcal N=D^{-\tfrac12}(D^2f):\Omega\longrightarrow
\{P\in{\rm Sym}^+(2):\det P=1\}\simeq\bbh^2,
\]
where ${\rm Sym}^+(2)$ denotes the set of real symmetric positive
definite $2\times2$ matrices.
Using this map and a bounded primitive of the hyperbolic area form,
we obtain an integral curvature estimate.
Combined with the positive function $D^{a+\tfrac12}$,
which is harmonic with respect to the Calabi metric $G$, this gives a spectral
inequality for the Gaussian curvature with coefficient greater
than $\tfrac14$ whenever $a\notin[-\tfrac23,-\tfrac13]$.
Under Calabi completeness, Castillon's theorem and the simple connectivity of $\Omega$
then imply that $(\Omega,G)$ is conformally equivalent to the plane.
The Liouville theorem then implies that $D$ is constant,
from which we deduce the following rigidity result.

\begin{thm}\label{1.0.1}
Let $f$ be a smooth strictly convex solution of the PDE (\ref{1.1}) on a convex domain in $\bbr^2$. If the
Calabi metric (\ref{1.3}) is complete and
\be\label{1.5}
a\notin\left[-\tfrac{2}{3},-\tfrac{1}{3}\right],
\ee
then $f$ is a quadratic polynomial.
\end{thm}

Let $\mathcal L_f(x)=(f_1(x),\ldots,f_n(x))^T$
denote the gradient map of $f$.  Throughout this paper, vectors are regarded as column vectors,
and the superscript $T$ denotes transpose.
The following theorem constructs non-quadratic Calabi complete
solutions on convex domains in every dimension $n\geq2$.

\begin{thm}\label{1.0.2}
Let $n\geq2$. The following functions are smooth strictly convex
non-quadratic solutions of (\ref{1.1}) on convex domains and have
complete Calabi metrics:
\begin{enumerate}
\item For $a\in[-\tfrac12,-\tfrac1{n+1}]$,
\[
f_a(x,y)=\int_1^x(x-t)(1+t)^{n-1}t^{\tfrac{1}{a}}\,dt
+\tfrac{|y|^2}{2(1+x)},
\qquad (x,y)\in(0,\infty)\times\bbr^{n-1}.
\]
\item For $b\in[-\tfrac n{n+1},-\tfrac12]$, set $a=-1-b$.
Let $u_b$ be the Legendre transform of $f_a$, defined on
$\mathcal L_{f_a}\bigl((0,\infty)\times\bbr^{n-1}\bigr)$.
It satisfies (\ref{1.1}) with exponent $b$.
Its explicit expression and convex domain are given in
(\ref{4.11}) and (\ref{4.12}), respectively.
\end{enumerate}
These two parameter intervals cover
$[-\tfrac n{n+1},-\tfrac1{n+1}]$.
\end{thm}

Together, Theorems \ref{1.0.1} and \ref{1.0.2} show that
non-quadratic Calabi complete solutions on convex domains
in dimension two exist precisely for
$a\in[-\tfrac23,-\tfrac13]$.
This establishes the sharpness of the rigidity range
in Theorem \ref{1.0.1}.

The paper is organized as follows.
Section 2 recalls the basic formulas of Calabi affine geometry
and establishes several preliminary lemmas.
Section 3 establishes a spectral inequality for the Gaussian
curvature and proves Theorem \ref{1.0.1}.
Section 4 constructs the non-quadratic solutions in
Theorem \ref{1.0.2} and verifies their Calabi completeness.

%%%%%%%%%%%%%%%%%%%%%%%
%%%%%%%%%%%%%%%%%%%%%%%%
\section{Preliminaries}

In this section, we recall the basic formulas of Calabi affine
geometry and establish some lemmas used below.
For background on Calabi affine geometry, see \cite{Ca1,LSZH}
and Section 3.3.4 of \cite{LXSJ}.

\vskip0.1in
\noindent{\textbf{2.1. Calabi affine geometry}}

Let $f$ be a strictly convex $C^\infty$-function on a convex domain
$\Omega\subset\bbr^n$, $n\geq2$.
Consider the graph hypersurface
\be\label{2.1}
M_f:=\{(x_1,\ldots,x_{n+1})\;|\;x_{n+1}=f(x_1,\ldots,x_n),\;(x_1,\ldots,x_n)\in\Omega\}.
\ee
Calabi affine geometry can also be essentially realised as a very special relative affine geometry of the graph hypersurface (\ref{2.1})
 by choosing the so-called {\em Calabi affine normalization} with $$Y=(0,\ldots,0,1)^T\in \bbr^{n+1}$$
being the fixed relative affine normal vector field.

For the graph immersion $\textbf{x}=(x_1,\ldots,x_n,f(x_1,\ldots,x_n))^T$, we have the decomposition
\be\label{2.2}
\tfrac{\partial^2\textbf{x}}{\partial x_i\partial x_j}
=\sum_k c_{ij}^k \textbf{x}_*\tfrac{\partial}{\partial x_k}+f_{ij}Y,
\ee
with respect to the bundle decomposition $\bbr^{n+1}=TM_f\oplus\bbr\cdot Y$, where the coefficients of the induced affine connection
satisfy $c_{ij}^k=0$. It follows that the relative affine metric is nothing but the \emph{Calabi metric}
\be\label{2.3} G=\sum \tfrac{\partial^2f}{\partial x_i\partial x_j}dx_idx_j.\ee
The Levi-Civita connection with respect to the metric $G$ has the Christoffel symbols
$$\hat \Gamma_{ij}^k=\tfrac{1}{2}\sum f^{kl}f_{ijl},$$
where, throughout the paper,
$$f_{ijk}=\tfrac{\partial^3f}{\partial x_i\partial x_j\partial x_k}, \quad (f^{kl})=(f_{ij})^{-1}. $$
Then we can rewrite the Gauss structure equation as follows:
\be\label{2.4} \textbf{x}_{,ij}=\sum A_{ij}^k\textbf{x}_*\tfrac{\partial }{\partial x_k}+G_{ij}Y,\ee
where the subscript $ _{,ij} $ denotes the covariant derivative with respect to
the Levi-Civita connection $\hat\nabla$ of the Calabi metric. The Fubini-Pick tensor $A_{ij}^k$ and the Weingarten tensor satisfy
\be\label{2.5} A_{ijk}=A_{ij}^lG_{kl}=-\tfrac{1}{2}f_{ijk},\quad B_{ij}=0,\ee
which means $A_{ijk}$ are symmetric in all indices.
The \emph{Tchebychev vector field} is
\be\label{2.6}
T:=\tfrac1n\sum_{i,j,l}G^{ij}A_{ij}^l\tfrac{\partial}{\partial x_l}
=-\tfrac1{2n}\nabla\ln D.
\ee
Here and throughout, gradients and tensor norms are taken
with respect to  the Calabi metric $G$ unless otherwise specified.
The Tchebychev operator is $\mathcal T:=\hat\nabla T$, and the \emph{relative Pick invariant} is
\be\label{2.7}
J:=\tfrac1{n(n-1)}\sum G^{il}G^{jp}G^{kq}A_{ijk}A_{lpq}.
\ee

We use $\Delta={\rm div}_G{\rm grad}_G$, so $-\Delta$ is nonnegative.
For a smooth function $\psi$, the Laplace-Beltrami operator
of $G$ takes the form
\be\label{2.8}
\Delta\psi
=\sum_{i,j}f^{ij}\psi_{ij}
-\tfrac12\sum_{i,j}f^{ij}(\ln D)_i\psi_j.
\ee

\vskip0.1in
\noindent{\textbf{2.2. The affine maximal type equation and related identities}}

We recall the following equivalent formulations of the affine
maximal type equation; see \cite{LXSJ} and Lemma 4.3 in \cite{SXX26}.

\begin{lem}[\cite{LXSJ,SXX26}]\label{2.0.1}
Let $M_f$ be a Calabi hypersurface with graph function $f$, and let $a\neq0$.
Then $f$ satisfies the PDE (\ref{1.1}) if and only if
\be\label{2.10}
{\rm tr}\mathcal{T}={\rm div}_G T=-\tfrac{1}{2n}\Delta\ln{D}=n(2a+1)|T|^2.
\ee
Equivalently, setting $c:=a+\tfrac12$, we obtain
\be\label{2.11}
\Delta\ln{D}=-4n^2c|T|^2=-c|\nabla\ln{D}|^2,\qquad \Delta (D^c)=0.
\ee
When $c=0$, the first identity reads $\Delta\ln D=0$,
while the second is trivial.
\end{lem}

The identities above yield the following integral formula.

\begin{lem}\label{2.0.2}
Let $f$ be a smooth strictly convex solution of (\ref{1.1})
on $\Omega$, and set $c:=a+\tfrac12$ and $v:=D^{c/2}>0$.
Then for every real-valued $\eta\in C_0^{\infty}(\Omega)$,
\be\label{2.12}
\int_\Omega|\nabla\eta|^2dV_G-\tfrac{c^2}{4}\int_\Omega|\nabla\ln{D}|^2\eta^2dV_G
=\int_\Omega v^2\left|\nabla\left(\tfrac{\eta}{v}\right)\right|^2dV_G\geq0.
\ee
\end{lem}

\textbf{Proof.} Since $v=D^{\tfrac{c}{2}}$, the chain rule and (\ref{2.11}) give
\be\label{2.13}
\tfrac{\Delta v}{v}
=\tfrac c2\Delta\ln D+\tfrac{c^2}{4}|\nabla\ln D|^2
=-\tfrac{c^2}{4}|\nabla\ln D|^2.
\ee
Set $\zeta:=\tfrac{\eta}{v}$. Then
$$|\nabla\eta|^2
=v^2|\nabla\zeta|^2+\lagl\nabla v,\nabla(v\zeta^2)\ragl.$$
Since $\eta$ has compact support and $v>0$ is smooth,
$v\zeta^2=\tfrac{\eta^2}{v}$ also has compact support. Integration by parts gives
$$\int_\Omega|\nabla\eta|^2dV_G
=\int_\Omega v^2|\nabla\zeta|^2dV_G
-\int_\Omega\tfrac{\Delta v}{v}\eta^2dV_G.$$
Substituting (\ref{2.13}) into this identity yields (\ref{2.12}). \hfill $\Box$

Let $f$ be a smooth strictly convex function on a convex domain
$\Omega\subset\bbr^n$.
Its Legendre transform $u$ is defined on
$\tilde\Omega:=\mathcal L_f(\Omega)$ by
\[
u(\xi)=x\cdot\xi-f(x),\qquad \xi=\mathcal L_f(x).
\]
For the dual potential $u$, write
\[
u_{ij}=\tfrac{\partial^2u}{\partial\xi_i\partial\xi_j},
\qquad \tilde D=\det(u_{ij}),
\qquad \tilde G=\sum_{i,j}u_{ij}d\xi_i d\xi_j.
\]
By Lemma 2.28 and the statement immediately following its proof
on p.~88 of \cite{LSZH},  the gradient map is a global
diffeomorphism onto its image.
By the Legendre duality formulas in Section 4 of \cite{Lian25},
with $\alpha=0$ for the metric, we obtain the isometry property
and the exponent transformation.
These facts are summarized in the following lemma.

\begin{lem}[\cite{Lian25,LSZH}]\label{2.0.3}
The map $\mathcal L_f:(\Omega,G)\to(\tilde\Omega,\tilde G)$ is a global isometry.
If $f$ satisfies (\ref{1.1}) with $a\neq0,-1$, then $u$ satisfies the same equation with exponent
\be\label{2.14}
\tilde a=-1-a.
\ee
In particular, $(\Omega,G)$ is complete if and only if
$(\tilde\Omega,\tilde G)$ is complete.
\end{lem}

To express the curvature in terms of the trace-free part of $A$, set
\be\label{2.15}
\stackrel{\circ}{A}{}_{ijk}:=A_{ijk}-\tfrac{n}{n+2}
\left(T_iG_{jk}+T_jG_{ik}+T_kG_{ij}\right).
\ee
Since $G^{jk}A_{ijk}=nT_i$, contraction gives $G^{jk}\stackrel{\circ}{A}_{ijk}=0$.
Taking norms and using $T=-\tfrac1{2n}\nabla\ln D$, we obtain
\be\label{2.16}
|A|^2=|\stackrel{\circ}{A}|^2+\tfrac{3n^2}{n+2}|T|^2
=|\stackrel{\circ}{A}|^2+\tfrac{3}{4(n+2)}|\nabla\ln D|^2.
\ee
We will also use the following curvature formulas.

\begin{lem}\label{2.0.4}
Let $M_f$ be a Calabi hypersurface of dimension $n\geq2$.
Then the scalar curvature $R$ of its Calabi metric satisfies
\be\label{2.17}
R=n(n-1)J-n^2|T|^2=|\stackrel{\circ}{A}|^2-\tfrac{n-1}{4(n+2)}|\nabla\ln{D}|^2.
\ee
In dimension two, the Gaussian curvature $K$ satisfies
\be\label{2.18}
K=\tfrac{1}{2}|\stackrel{\circ}{A}|^2-\tfrac{1}{32}|\nabla\ln{D}|^2.
\ee
\end{lem}

The first equality in (\ref{2.17}) is the standard scalar
curvature identity; see Section 3.3.4 in \cite{LXSJ}.
The second follows from (\ref{2.16}), while (\ref{2.18})
follows from $R=2K$.

\vskip0.1in
\noindent{\textbf{2.3. Theorem of Castillon}}

In this subsection, we recall Theorem A of Castillon \cite{Cast}.
Let $(M,G)$ be a connected, oriented, complete noncompact Riemannian
surface with Gaussian curvature $K$. For $\lambda\in\bbr$, define
\be\label{3.36}
q_\lambda(\eta)
:=\int_M\left(|\nabla\eta|^2+\lambda K\eta^2\right)dV_G,
\qquad \eta\in C_0^\infty(M).
\ee
Integration by parts identifies $q_\lambda$ as the quadratic
form associated with $-\Delta+\lambda K$:
$$
q_\lambda(\eta)=\int_M\eta(-\Delta+\lambda K)\eta\,dV_G.
$$
Castillon uses the nonnegative Laplacian
$\Delta_{\rm Cast}=-\Delta$; see pp.~272--273 in \cite{Cast}.
Thus his operator $\Delta_{\rm Cast}+\lambda K$ agrees with ours.
Define
\be\label{3.37}
b_G:=\sup\{\lambda\in\bbr\;|\; q_\lambda(\eta)\geq0\ \hbox{ for all }\eta\in C_0^{\infty}(M)\}.
\ee
Castillon's theorem can then be stated as follows.

\begin{thm}[\cite{Cast}]\label{3.0.4}
Let $(M^2,G)$ be a connected, oriented, complete noncompact Riemannian surface with 
Gaussian curvature $K$, and let $b_G$ be defined by (\ref{3.37}). 
If $b_G>\tfrac14$, then $M$ is conformally equivalent to $\bbc$ or $\bbc\setminus\{0\}$.
\end{thm}

%%%%%%%%%%%%%%%%%%%%%%%
%%%%%%%%%%%%%%%%%%%%%%%%
\section{Sharp Bernstein rigidity in dimension two}
\label{3.0}

We prove Theorem \ref{1.0.1} using an integral curvature estimate
derived from the normalized Hessian map.
Combined with Lemma \ref{2.0.2}, this gives a spectral
inequality to which Castillon's theorem applies.
A Liouville argument then yields rigidity.
Throughout this section, we set $n=2$, let the indices $i,j,k$
take values in $\{1,2\}$, and orient $\Omega$ by $dx_1\wedge dx_2$.

\vskip0.1in
\noindent{\textbf{3.1. The normalized Hessian map and curvature estimates}}

We use a matrix model of the hyperbolic plane to study the
normalized Hessian map; see Remark 3.2(b) in \cite{DP}.
Let ${\rm Sym}(2)$ denote the vector space of real symmetric $2\times2$ matrices, and let
${\rm Sym}^+(2)$ be the subset of positive definite matrices.
Set
$$\mathcal P:=\{P\in{\rm Sym}^+(2)\mid\det P=1\}.$$
We identify ${\rm Sym}(2)$ with $\bbr^3$ via
$$
P=\lmx p_{11}&p_{12}\\p_{12}&p_{22}\rmx
\longmapsto(p_{11},p_{12},p_{22})^T.
$$
Then $P\in\mathcal P$ if and only if
\be\label{3.01}
p_{11}>0,\qquad p_{11}p_{22}-p_{12}^2=1.
\ee
The second condition in (\ref{3.01}) gives $p_{22}=\tfrac{1+p_{12}^2}{p_{11}}$.
Thus $\mathcal P$ is parametrized by
$$
\Phi:(0,\infty)\times\bbr\longrightarrow\mathcal P,\qquad
(p_{11},p_{12})\longmapsto
\lmx p_{11}&p_{12}\\p_{12}&\tfrac{1+p_{12}^2}{p_{11}}\rmx.
$$
In particular, $\mathcal P$ is the graph of the smooth function
$p_{22}=\tfrac{1+p_{12}^2}{p_{11}}$ on $(0,\infty)\times\bbr$,
and hence a smooth embedded surface with global coordinates
$(p_{11},p_{12})$.

To determine $T_P\mathcal P$, let $V\in T_P\mathcal P$ and
choose a smooth curve $P(t)$ in $\mathcal P$ with
$P(0)=P$ and $P'(0)=V$.
Differentiating $\det P(t)=1$ gives
$$
0=\left.\tfrac{d}{dt}\det P(t)\right|_{t=0}
=\det P\,{\rm tr}(P^{-1}V)
={\rm tr}(P^{-1}V),
$$
where ${\rm tr}$ denotes the matrix trace.
The condition ${\rm tr}(P^{-1}V)=0$ defines a two-dimensional
subspace of ${\rm Sym}(2)$.
Since $T_P\mathcal P$ is contained in this subspace and
has the same dimension, we obtain
\be\label{3.1}
T_P\mathcal P
=\{V\in{\rm Sym}(2)\mid{\rm tr}(P^{-1}V)=0\}.
\ee

For $V,W\in T_P\mathcal P$, define
\be\label{3.2}
\lagl V,W\ragl_P=\tfrac12{\rm tr}(P^{-1}VP^{-1}W).
\ee
This defines a smooth symmetric bilinear form on $\mathcal P$.
We verify positive definiteness by computing the coefficients
in coordinates. This also identifies the resulting metric
with the standard hyperbolic metric.
Introduce new coordinates
\be\label{3.3'}
X=\tfrac{p_{12}}{p_{11}},\qquad Y=\tfrac1{p_{11}}>0.
\ee
Then $p_{11}=\tfrac{1}{Y}$, $p_{12}=\tfrac{X}{Y}$, and the determinant condition in (\ref{3.01}) gives
$p_{22}=\tfrac{X^2+Y^2}{Y}$.
Thus the parametrization becomes
\be\label{3.3}
P(X,Y)=\tfrac1Y\lmx 1&X\\X&X^2+Y^2\rmx,
\qquad (X,Y)\in\bbr\times(0,\infty).
\ee
The matrices $P_X:=\tfrac{\partial P}{\partial X}$ and
$P_Y:=\tfrac{\partial P}{\partial Y}$ form a basis of
$T_P\mathcal P$.
A direct calculation gives
\be\label{3.4}P^{-1}=\tfrac1Y\lmx X^2+Y^2&-X\\-X&1\rmx,\;
P_X=\tfrac1Y\lmx0&1\\1&2X\rmx,\;
P_Y=\tfrac1{Y^2}\lmx-1&-X\\-X&Y^2-X^2\rmx.\ee
Substituting (\ref{3.4}) into (\ref{3.2}), we obtain
$$
\lagl P_X,P_X\ragl_P
=\lagl P_Y,P_Y\ragl_P=Y^{-2},\qquad
\lagl P_X,P_Y\ragl_P=0.
$$
Thus (\ref{3.2}) defines a Riemannian metric $G_{\mathcal P}$, which in the global coordinates $(X,Y)$ takes the form
\be\label{3.5}
G_{\mathcal P}=\tfrac{dX^2+dY^2}{Y^2}.
\ee
This is the standard upper half-plane metric of curvature $-1$,
so $(\mathcal P,G_{\mathcal P})$ is isometric to the hyperbolic
plane $\bbh^2$.

We give $\mathcal P$ the orientation determined by $dX\wedge dY$. Its area form $\omega$
has a bounded primitive $\alpha$, given by
\be\label{3.6}
\omega=\tfrac{dX\wedge dY}{Y^2},\qquad
\alpha=\tfrac{dX}{Y},\qquad
d\alpha=\omega,\qquad |\alpha|_{G_{\mathcal P}}=1.
\ee
Since $\det(D^{-\tfrac12}D^2f)=1$, we define the
normalized Hessian map by
\be\label{3.7}
\mathcal N:\Omega\longrightarrow\mathcal P,\qquad
x\longmapsto D^{-\tfrac12}(D^2f)(x).
\ee
The following proposition gives formulas for the signed Jacobian
and the operator norm of $d\mathcal N$.

\begin{prop}\label{3.0.1}
Let $f$ be a smooth strictly convex function on a domain
$\Omega\subset\bbr^2$, and let $G$ be its Calabi metric.
Then the normalized Hessian map $\mathcal N$ defined in
(\ref{3.7}) satisfies
\be\label{3.8}
\mathcal{N}^*\omega=-2K\,dV_G,\qquad
\|d\mathcal{N}\|=\tfrac{1}{4}|\nabla\ln{D}|+|\stackrel{\circ}{A}|,
\ee
where $K$ is the Gaussian curvature of $G$ and $\|d\mathcal N\|=\sup_{|v|_G=1}|d\mathcal N(v)|_{G_{\mathcal P}}$ is the operator norm of $d\mathcal N$.
\end{prop}

\textbf{Proof.}
We first verify that (\ref{3.8}) is invariant under affine coordinate changes that preserve orientation.
Consider the affine coordinate change
\[
x=\chi(z):=x_0+Lz,
\]
where $x_0\in\bbr^2$ and $L\in GL(2,\bbr)$ are fixed,
with $\det L>0$.
On $\bar\Omega:=\chi^{-1}(\Omega)$, define
\[
\bar f:=f\circ\chi,\qquad
\bar D:=\det(D_z^2\bar f),\qquad
\bar{\mathcal N}:=\bar D^{-\tfrac12}(D_z^2\bar f).
\]
Applying the chain rule and taking determinants, we obtain
\be\label{3.9}
D_z^2\bar f=L^T(D_x^2f)L,\qquad
\bar D=(\det L)^2D,\qquad
\bar{\mathcal N}=M^T\mathcal N M,
\ee
where $M:=(\det L)^{-\tfrac12}L\in SL(2,\bbr)$.
The subscripts $x$ and $z$ indicate the differentiation variables.
Here and below, unbarred quantities are evaluated at $x=\chi(z)$
unless otherwise indicated.

Let $\bar G$ be the Calabi metric of $\bar f$.
The map $\chi:(\bar\Omega,\bar G)\to(\Omega,G)$
extends to an affine transformation of $\bbr^3$ that fixes
the Calabi affine normal and maps the two graphs by
\[
(z,\bar f(z))\longmapsto
(x_0+Lz,\bar f(z))=(x,f(x)).
\]
The invariance of the Calabi metric and the Fubini-Pick form gives
\be\label{3.10}
\bar G=\chi^*G,\qquad \bar A=\chi^*A.
\ee
Thus $\chi$ is an isometry.
Since $\det L>0$, it also preserves orientation. Hence
\be\label{3.11}
\bar K=K\circ\chi,\qquad dV_{\bar G}=\chi^*dV_G.
\ee
Taking trace-free parts in (\ref{3.10}) gives
\be\label{3.12}
\stackrel{\circ}{\bar A}=\chi^*\stackrel{\circ}{A},
\qquad
|\stackrel{\circ}{\bar A}|_{\bar G}
=|\stackrel{\circ}{A}|_G\circ\chi.
\ee
Since $\ln\bar D=(\ln D)\circ\chi+2\ln\det L$ and
$\det L$ is constant, the isometry $\chi$ gives
\be\label{3.13}
|\nabla_{\bar G}\ln\bar D|_{\bar G}
=|\nabla_G\ln D|_G\circ\chi.
\ee

For $P\in\mathcal P$, the matrix $M^TPM$ is positive definite
and $\det(M^TPM)=1$.
We next show that the map
$$
\Psi_M:\mathcal P\longrightarrow\mathcal P,\qquad
\Psi_M(P)=M^TPM,
$$
preserves the metric $G_{\mathcal P}$.
For $V\in T_P\mathcal P$, choose a smooth curve $P(t)$ in
$\mathcal P$ with $P(0)=P$ and $P'(0)=V$.
Since $M$ is fixed,
$$
d(\Psi_M)_P(V)
=\left.\tfrac{d}{dt}\bigl(M^TP(t)M\bigr)\right|_{t=0}
=M^TVM.
$$
For $V,W\in T_P\mathcal P$, using
$(M^TPM)^{-1}=M^{-1}P^{-1}(M^T)^{-1}$ in (\ref{3.2}) yields
\be\label{3.14}
\begin{aligned}
\lagl d(\Psi_M)_P(V),d(\Psi_M)_P(W)\ragl_{\Psi_M(P)}
&=\tfrac12{\rm tr}\bigl(M^{-1}P^{-1}VP^{-1}WM\bigr)\\
&=\tfrac12{\rm tr}(P^{-1}VP^{-1}W)\\
&=\lagl V,W\ragl_P.
\end{aligned}
\ee
The second equality follows from cyclicity of the trace.
Since $\Psi_M$ has the smooth inverse $\Psi_{M^{-1}}$,
it is an isometry of $(\mathcal P,G_{\mathcal P})$.

The map $\Psi_M$ also preserves orientation.
Indeed, since $SL(2,\bbr)$ is path connected, there is a
continuous path $M_t\in SL(2,\bbr)$, $0\leq t\leq1$,
with $M_0=I$ and $M_1=M$.
Since each $\Psi_{M_t}$ is a diffeomorphism, its signed
Jacobian is nonzero.
At each fixed $P\in\mathcal P$, this Jacobian depends
continuously on $t$.
It is positive at $t=0$, where $\Psi_{M_0}$ is the identity,
and is therefore positive for all $t\in[0,1]$.
Hence $\Psi_M$ preserves orientation and
\be\label{3.15}
\Psi_M^*\omega=\omega.
\ee

By (\ref{3.9}), $\bar{\mathcal N}=\Psi_M\circ\mathcal N\circ\chi$.
Using (\ref{3.15}), we obtain
\be\label{3.16}
\bar{\mathcal N}^*\omega
=\chi^*\bigl(\mathcal N^*(\Psi_M^*\omega)\bigr)
=\chi^*(\mathcal N^*\omega).
\ee
By (\ref{3.10}) and (\ref{3.14}), the chain rule gives
\be\label{3.17}
\|d\bar{\mathcal N}_z\|
=\sup_{|v|_{\bar G}=1}
|d\mathcal N_{\chi(z)}(d\chi_z(v))|_{G_{\mathcal P}}
=\|d\mathcal N_{\chi(z)}\|.
\ee
Equations (\ref{3.11}) and (\ref{3.16})
show that the first identity in (\ref{3.8}) is invariant under this coordinate change.
Equations (\ref{3.12}), (\ref{3.13})
and (\ref{3.17}) give the same conclusion for the second.

Fix a point $p\in\Omega$ and choose positively oriented affine
coordinates such that
\be\label{3.18}
G_{ij}(p)=\delta_{ij},\quad
(\ln D)_1(p)=\sigma,\quad
(\ln D)_2(p)=0,\quad
\sigma=|\nabla\ln D|(p)\geq0.
\ee
At $p$, the trace-freeness and symmetry of
$\stackrel{\circ}{A}$ imply
$$
\stackrel{\circ}{A}_{111}+\stackrel{\circ}{A}_{122}=0,\qquad
\stackrel{\circ}{A}_{112}+\stackrel{\circ}{A}_{222}=0.
$$
Setting $b_1=\stackrel{\circ}{A}_{111}$ and
$b_2=\stackrel{\circ}{A}_{112}$, we obtain
\be\label{3.19}
\stackrel{\circ}{A}_{111}=-\stackrel{\circ}{A}_{122}=b_1,\qquad
\stackrel{\circ}{A}_{112}=-\stackrel{\circ}{A}_{222}=b_2,\qquad
|\stackrel{\circ}{A}|^2=4(b_1^2+b_2^2).
\ee
Since $T_1=-\tfrac{\sigma}{4}$ and $T_2=0$ at $p$,
(\ref{2.15}) and $f_{ijk}=-2A_{ijk}$ give
\be\label{3.20}
(f_{jk1})=\lmx
\tfrac{3\sigma}{4}-2b_1&-2b_2\\
-2b_2&\tfrac{\sigma}{4}+2b_1
\rmx,\qquad
(f_{jk2})=\lmx
-2b_2&\tfrac{\sigma}{4}+2b_1\\
\tfrac{\sigma}{4}+2b_1&2b_2
\rmx.
\ee
Differentiating (\ref{3.7}) and using $(D^2f)(p)=I$ and $D(p)=1$,
we obtain at $p$
\be\label{3.21}
\tfrac{\partial\mathcal N}{\partial x_1}
=(f_{ij1})-\tfrac \sigma2 I,\qquad
\tfrac{\partial\mathcal N}{\partial x_2}
=(f_{ij2}).
\ee

We compute $\mathcal N^*\omega$ using the coordinates $(X,Y)$
on $\mathcal P$. Writing $\mathcal N=(\mathcal N_{ij})$, (\ref{3.3'}) gives
$$
X\circ\mathcal N=\tfrac{\mathcal N_{12}}{\mathcal N_{11}},
\qquad
Y\circ\mathcal N=\tfrac1{\mathcal N_{11}}.
$$
Differentiating the preceding identities and using
$\mathcal N(p)=I$, (\ref{3.20}) and (\ref{3.21}),
we obtain at $p$
$$
\begin{aligned}
d(X\circ\mathcal N)
&=-2b_2\,dx_1+\left(\tfrac{\sigma}{4}+2b_1\right)dx_2,\\
d(Y\circ\mathcal N)
&=-\left(\tfrac{\sigma}{4}-2b_1\right)dx_1+2b_2\,dx_2.
\end{aligned}
$$
Using $\omega=\tfrac1{Y^2}dX\wedge dY$ and
$(Y\circ\mathcal N)(p)=1$, together with
$G_{ij}(p)=\delta_{ij}$ and (\ref{2.18}), we obtain
\be\label{3.24}
\begin{aligned}
(\mathcal N^*\omega)_p
&=\left[\tfrac{1}{(Y\circ\mathcal N)^{2}}
d(X\circ\mathcal N)\wedge d(Y\circ\mathcal N)\right]_p\\
&=\left[\tfrac{{\sigma}^2}{16}-4(b_1^2+b_2^2)\right]
(dx_1\wedge dx_2)_p\\
&=\left(\tfrac1{16}|\nabla\ln D|^2
-|\stackrel{\circ}{A}|^2\right)(p)\,(dV_G)_p\\
&=-2K(p)(dV_G)_p.
\end{aligned}
\ee
Since $p$ is arbitrary, this proves the first identity in (\ref{3.8}).

We next compute the operator norm of $d\mathcal N_p$.
Since $\mathcal N(p)=I$, the differential of $\mathcal N$ at $p$ is
\[
d\mathcal N_p:T_p\Omega\longrightarrow T_I\mathcal P.
\]
Setting $P=I$ in (\ref{3.1}), we obtain
\[
T_I\mathcal P
=\{V\in{\rm Sym}(2):{\rm tr}(V)=0\}.
\]
We identify $T_I\mathcal P$ with $\bbr^2$ by
\be\label{3.22}
V=\lmx u&v\\v&-u\rmx\longmapsto(u,v)^T.
\ee
By (\ref{3.2}),
$$
|V|_I^2=\tfrac12{\rm tr}(V^2)=u^2+v^2.
$$
Thus (\ref{3.22}) is an isometry from
$(T_I\mathcal P,\langle\cdot,\cdot\rangle_I)$
to $\bbr^2$ with its Euclidean inner product.
Under this identification, (\ref{3.21}) gives
the matrix of $d\mathcal N_p$ as
\be\label{3.23}
(d\mathcal N_p)
=\lmx \tfrac{\sigma}{4}-2b_1&-2b_2\\
-2b_2&\tfrac{\sigma}{4}+2b_1\rmx.
\ee
The symmetric matrix (\ref{3.23}) has eigenvalues
$\tfrac{\sigma}{4}\pm2\sqrt{b_1^2+b_2^2}$.
By the spectral theorem, its operator norm is the larger absolute eigenvalue.
Since $G_{ij}(p)=\delta_{ij}$, (\ref{3.22}) is an isometry,
and $\sigma\geq0$, we obtain
$$
\|d\mathcal N_p\|
=\tfrac{\sigma}{4}+2\sqrt{b_1^2+b_2^2}
=\tfrac14|\nabla\ln D|(p)+|\stackrel{\circ}{A}|(p).
$$
As $p$ is arbitrary, this proves the second identity in (\ref{3.8}).
\hfill $\Box$

Using the bounded primitive $\alpha$ in (\ref{3.6}), we derive
the following integral curvature estimate from Proposition \ref{3.0.1}.

\begin{lem}\label{3.0.2}
Let $\eta\in C_0^\infty(\Omega)$ be real-valued, and define
\be\label{3.25}
E:=\int_\Omega|\nabla\eta|^2dV_G,\qquad
P:=\int_\Omega|\nabla\ln D|^2\eta^2dV_G,
\ee
\be\label{3.26}
Q:=\int_\Omega|\stackrel{\circ}{A}|^2\eta^2dV_G,\qquad
\mathcal K:=\int_\Omega K\eta^2dV_G.
\ee
Then
\be\label{3.27}
\mathcal K=\tfrac{Q}{2}-\tfrac{P}{32},\qquad
|\mathcal K|\leq\sqrt E\left(\tfrac{\sqrt P}{4}+\sqrt Q\right).
\ee
\end{lem}

\textbf{Proof.}
Multiplying (\ref{2.18}) by $\eta^2$ and integrating over $\Omega$
gives the first identity in (\ref{3.27}).
To prove the inequality, set $\theta:=\mathcal N^*\alpha$,
where $\alpha$ is the bounded primitive in (\ref{3.6}).
Since exterior differentiation commutes with pullback,
(\ref{3.6}) and (\ref{3.8}) give
\be\label{3.28}
d\theta=d(\mathcal N^*\alpha)=\mathcal N^*(d\alpha)=\mathcal N^*\omega=-2K\,dV_G.
\ee
Since $|\alpha|_{G_{\mathcal P}}=1$ by (\ref{3.6}),
for any tangent vector $v$ with $|v|_G=1$, we have
$$
\begin{aligned}
|\theta(v)|
=|\alpha(d\mathcal N(v))|\leq|\alpha|_{G_{\mathcal P}}\,|d\mathcal N(v)|_{G_{\mathcal P}}
\leq\|d\mathcal N\|.
\end{aligned}
$$
Taking the supremum over $|v|_G=1$ and using (\ref{3.8}) gives
\be\label{3.29}
|\theta|_G
=\sup_{|v|_G=1}|\theta(v)|
\leq\|d\mathcal N\|
=\tfrac14|\nabla\ln D|+|\stackrel{\circ}{A}|.
\ee
The one-form $\eta^2\theta$ has compact support. Stokes' theorem and
(\ref{3.28}) therefore imply
$$0=\int_\Omega d(\eta^2\theta)
=2\int_\Omega\eta\,d\eta\wedge\theta-2\int_\Omega K\eta^2dV_G.$$
Consequently,
\be\label{3.30}
\mathcal K=\int_\Omega\eta\,d\eta\wedge\theta.
\ee
Using $|d\eta\wedge\theta|_G\leq|d\eta|_G|\theta|_G$
and $|d\eta|_G=|\nabla\eta|$, we obtain from
(\ref{3.29}) and (\ref{3.30}) that
$$
|\mathcal K|
\leq\int_\Omega|\eta|\,|\nabla\eta|\,|\theta|_G\,dV_G
\leq\int_\Omega|\eta|\,|\nabla\eta|
\left(\tfrac14|\nabla\ln D|+|\stackrel{\circ}{A}|\right)dV_G.
$$
Applying the Cauchy--Schwarz inequality to each term gives
$$|\mathcal K|
\leq\tfrac14\sqrt E\sqrt P+\sqrt E\sqrt Q
=\sqrt E\left(\tfrac{\sqrt P}{4}+\sqrt Q\right),$$
which proves the inequality in (\ref{3.27}). \hfill $\Box$

\vskip0.1in
\noindent{\textbf{3.2. A spectral inequality for the Gaussian curvature}}

Combining Lemmas \ref{2.0.2} and \ref{3.0.2},
we obtain the following spectral inequality for the Gaussian curvature.

\begin{prop}\label{3.0.3}
Let $f$ be a smooth strictly convex solution of (\ref{1.1})
on a domain $\Omega\subset\bbr^2$, and suppose that
$c:=a+\tfrac12\neq0$. Set
\be\label{3.31}
\delta:=\min\left\{|c|,\tfrac14\right\},\qquad
\lambda:=\tfrac{\delta}{1-2\delta}.
\ee
Then, for every real-valued $\eta\in C_0^\infty(\Omega)$,
\be\label{3.32}
\int_\Omega\left(|\nabla\eta|^2+\lambda K\eta^2\right)dV_G\geq0.
\ee
In particular, if $a\notin[-\tfrac23,-\tfrac13]$,
then (\ref{3.32}) holds with $\lambda>\tfrac14$.
\end{prop}

\textbf{Proof.} We use the notation in (\ref{3.25}) and (\ref{3.26}).
Since $0<\delta\leq\tfrac14$, we have $\lambda>0$.
As $E\geq0$, (\ref{3.32}) holds whenever
$\mathcal K\geq0$. It therefore suffices to consider
the case $\mathcal K<0$.
Set $s:=\sqrt E$, $p:=\sqrt P$ and $q:=\sqrt Q$.
Since $\mathcal K<0$, the first identity in (\ref{3.27}) gives
$p>4q\geq0$.
Using (\ref{2.12}) and $\delta\leq|c|$, we obtain
\be\label{3.33}
s^2\geq\tfrac{c^2}{4}p^2\geq\tfrac{\delta^2}{4}p^2>0.
\ee
Moreover, (\ref{3.27}) gives
$$
\tfrac{(p-4q)(p+4q)}{32}
=-\mathcal K
\leq\tfrac{s(p+4q)}{4}.
$$
Dividing by $p+4q>0$ and using (\ref{3.33}), we obtain
\be\label{3.34}
p\leq\tfrac{2s}{\delta},\qquad p-4q\leq8s.
\ee
We claim that
\be\label{3.35}
-\mathcal{K}\leq\left(\tfrac{1}{\delta}-2\right)E.
\ee
If $p\leq8s$, then
$$
-\mathcal K=\tfrac{p^2-16q^2}{32}
\leq2s^2
\leq\left(\tfrac1\delta-2\right)E,
$$
where the last inequality follows from $\delta\leq\tfrac14$ and $s^2=E$.
If $p>8s$, the second inequality in (\ref{3.34}) gives
$4q\geq p-8s>0$, hence $16q^2\geq(p-8s)^2$. Therefore,
$$
-\mathcal K
=\tfrac{p^2-16q^2}{32}
\leq\tfrac{p^2-(p-8s)^2}{32}
=\tfrac{ps}{2}-2s^2
\leq\left(\tfrac1\delta-2\right)E,
$$
where the last inequality follows from the first inequality in
(\ref{3.34}) and $s^2=E$.
Thus (\ref{3.35}) holds in both cases.
Since $\lambda=(\tfrac1\delta-2)^{-1}>0$, multiplying
(\ref{3.35}) by $\lambda$ gives $E+\lambda\mathcal K\geq0$, which proves (\ref{3.32}).
Finally, when $a\notin[-\tfrac23,-\tfrac13]$, we have
$|c|>\tfrac16$.
The definition of $\delta$ then gives
$\delta>\tfrac16$,
and hence $\lambda>\tfrac14$.
\hfill $\Box$

\vskip0.1in
\noindent{\textbf{3.3. Proof of Theorem \ref{1.0.1}}}

Recall the definition of$$
q_\lambda(\eta)
:=\int_M\left(|\nabla\eta|^2+\lambda K\eta^2\right)dV_G,
\qquad \eta\in C_0^\infty(M) 
$$
and 
$$
b_G:=\sup\{\lambda\in\bbr\;|\; q_\lambda(\eta)\geq0\ \hbox{ for all }\eta\in C_0^{\infty}(M)\}.
$$
By Proposition \ref{3.0.3}, since
$a\notin[-\tfrac23,-\tfrac13]$, there exists $\lambda>\tfrac14$
such that $q_\lambda(\eta)\geq0$ for every
$\eta\in C_0^\infty(\Omega)$.
Hence $b_G\geq\lambda>\tfrac14$ by the definition of $b_G$.
Theorem \ref{3.0.4} therefore implies that $(\Omega,G)$ is
conformally equivalent to $\bbc$ or $\bbc\setminus\{0\}$.
Since $\Omega$ is convex, it is simply connected and therefore
cannot be conformally equivalent to $\bbc\setminus\{0\}$.
Set $c:=a+\tfrac12\neq0$.
By Lemma \ref{2.0.1}, $U:=D^c$ is positive and harmonic
on $(\Omega,G)$.
Let $\mathcal F:\bbc\to(\Omega,G)$ be a conformal diffeomorphism.
By conformal invariance of harmonicity in dimension two,
$U\circ\mathcal F$ is a positive harmonic function on $\bbc$.
The classical Liouville theorem implies that $U\circ\mathcal F$,
and hence $U$, is constant. Since $c\neq0$, $D$ is also constant.
Together with the completeness of $G$, the first corollary to Theorem 1 of \cite{Ca1}, 
implies that $f$ is a quadratic polynomial.
This completes the proof of Theorem \ref{1.0.1}.
\hfill $\Box$

%%%%%%%%%%%%%%%%%%%%%%%
%%%%%%%%%%%%%%%%%%%%%%%%
\section{Calabi complete non-quadratic examples in dimensions $n\geq2$}
\label{4.0}

We first construct explicit solutions and determine the range
of parameters for which they are Calabi complete.
Their Legendre transforms provide the remaining examples
in Theorem \ref{1.0.2}.

\vskip0.1in
\noindent{\textbf{4.1. An explicit family of solutions and its completeness range}}

Write $x=x_1$ and $y=(x_2,\ldots,x_n)^T$. For a real parameter $a\neq0$, define
$f_a$ on $\Omega:=(0,\infty)\times\bbr^{n-1}$ by
\be\label{4.1}
f_a(x,y)=\vfi_a(x)+\tfrac{|y|^2}{2(1+x)},\qquad
\vfi_a(x)=\int_1^x(x-t)(1+t)^{n-1}t^{\tfrac{1}{a}}\,dt.
\ee
The integral is well defined for every $x>0$, with
$\int_1^x=-\int_x^1$ when $x<1$.
Setting $H_a:=\vfi_a'$, differentiation gives
\be\label{4.2}
H_a(x)=\int_1^x(1+t)^{n-1}t^{\tfrac{1}{a}}\,dt,\qquad
H_a'(x)=(1+x)^{n-1}x^{\tfrac{1}{a}}>0.
\ee

\begin{prop}\label{4.0.1}
For every $a\neq0$, the function $f_a$ defined by (\ref{4.1})
is a smooth strictly convex non-quadratic solution of (\ref{1.1})
on $(0,\infty)\times\bbr^{n-1}$ with exponent $a$.
Its Calabi metric is complete if and only if
\be\label{4.3}
a\in\left[-\tfrac12,-\tfrac1{n+1}\right].
\ee
\end{prop}

\textbf{Proof.} Write $f=f_a$. Put $r:=1+x$ and $\mu(x):=H_a'(x)=r^{n-1}x^{\tfrac{1}{a}}$.
Differentiating (\ref{4.1}) gives
$$
f_x=H_a(x)-\tfrac{|y|^2}{2r^2},\qquad
f_{y_j}=\tfrac{y_j}{r},
$$
and hence
\be\label{4.4}
(D^2f)=\lmx
\mu+\tfrac{|y|^2}{r^3}&-\tfrac{y^T}{r^2}\\
-\tfrac{y}{r^2}&\tfrac1r I_{n-1}
\rmx.
\ee
Here $I_{n-1}$ denotes
the $(n-1)\times(n-1)$ identity matrix.
For an arbitrary vector $w=\binom{s}{v}\in\bbr^n$, with
$s\in\bbr$ and $v\in\bbr^{n-1}$, (\ref{4.4}) gives
\be\label{4.5}
w^T(D^2f)w
=\mu s^2+\tfrac1r\left|v-\tfrac{y}{r}s\right|^2.
\ee
Since $\mu>0$ and $r>0$, the right-hand side is positive for
$w\neq0$, so $f$ is strictly convex.
Taking the determinant in (\ref{4.4}), we obtain
\be\label{4.6}
D=\det(D^2f)=r^{-(n-1)}\mu=x^{\tfrac{1}{a}}.
\ee
Since $D^a=x$, $f$ satisfies (\ref{1.1}).
Moreover, $D$ is nonconstant, so $f$ is non-quadratic.

To study completeness, introduce the global coordinate change
$$
(x,y)\longmapsto(x,z),\qquad z=\tfrac{y}{1+x}.
$$
Since $r=1+x$, we obtain $dy=z\,dx+r\,dz$.
Thus (\ref{4.5}) gives
\be\label{4.7}
G=\mu(x)\,dx^2+\tfrac1r\left|dy-\tfrac yr\,dx\right|^2=\mu(x)\,dx^2+r|dz|^2.
\ee
We claim that $G$ is complete if and only if both integrals
\be\label{4.8}
\int_0^1\sqrt{\mu(x)}\,dx,\qquad
\int_1^\infty\sqrt{\mu(x)}\,dx
\ee
diverge.
Necessity is immediate by considering curves with constant $z$.
To prove sufficiency, suppose that both integrals in (\ref{4.8}) diverge.
Define
$$
t(x):=\int_1^x\sqrt{\mu(s)}\,ds.
$$
Since $t'(x)=\sqrt{\mu(x)}>0$ and
$$
\lim_{x\to0^+}t(x)=-\infty,\qquad
\lim_{x\to\infty}t(x)=+\infty,
$$
$t$ is a diffeomorphism from $(0,\infty)$ onto $\bbr$.
Writing its inverse as $x=x(t)$, we obtain global coordinates
$(t,z)\in\bbr^n$.
As $dt=\sqrt{\mu(x)}\,dx$, (\ref{4.7}) becomes
\be\label{4.9}
G=dt^2+(1+x(t))|dz|^2
\geq dt^2+|dz|^2=:g_E.
\ee
Since $g_E$ is the complete Euclidean metric on $\bbr^n$,
$G$ is complete as well.

To determine the range of $a$ for which $G$ is complete, note that
$$
\sqrt{\mu(x)}=(1+x)^{\tfrac{n-1}{2}}x^{\tfrac{1}{2a}}.
$$
Hence
$$
\sqrt{\mu(x)}\sim x^{\tfrac{1}{2a}}\quad(x\to0^+),\qquad
\sqrt{\mu(x)}\sim x^{\tfrac{n-1}{2}+\tfrac{1}{2a}}\quad(x\to\infty).
$$
Thus the first integral in (\ref{4.8}) diverges precisely when
$\tfrac{1}{2a}\leq-1$, that is, $-\tfrac12\leq a<0$.
Under this necessary condition, the second integral diverges precisely when
$\tfrac{n-1}{2}+\tfrac{1}{2a}\geq-1$, that is, $a\leq-\tfrac1{n+1}$.
Thus both integrals diverge precisely for
$a\in[-\tfrac12,-\tfrac1{n+1}]$, proving (\ref{4.3}).
\hfill $\Box$

\vskip0.1in
\noindent{\textbf{4.2. Legendre dual solutions and proof of Theorem \ref{1.0.2}}}

For a smooth strictly convex function on a convex domain,
the gradient map is a diffeomorphism onto its image, and
the Legendre transform has a positive definite Hessian;
see Section 2.7.1, pp.~88--89, of \cite{LSZH}.
Note that the image of the gradient map need not be convex.
We therefore also verify the convexity of the dual domain
to obtain a strictly convex solution on a convex domain,
as required by Theorem \ref{1.0.2}.

Let $a\in[-\tfrac12,-\tfrac1{n+1}]$ and set $b=-1-a$.
We compute the Legendre transform $u_b$ of $f_a$ in (\ref{4.1}).
Writing $\xi=\binom{\xi_1}{\xi'}\in\bbr\times\bbr^{n-1}$,
the gradient map $\xi=\mathcal L_{f_a}(x,y)$ is given by
\be\label{4.10}
\xi_1=H_a(x)-\tfrac{|y|^2}{2(1+x)^2},
\qquad
\xi'=\tfrac{y}{1+x}.
\ee
Consequently,
\be\label{4.10'}
\xi_1+\tfrac{|\xi'|^2}{2}=H_a(x),
\qquad
y=(1+x)\xi'.
\ee
On the gradient image
$\tilde\Omega_a:=\mathcal L_{f_a}(\Omega)$,
the Legendre transform is given by
\be\label{4.11}
\begin{aligned}
u_b(\xi)
&=x\xi_1+y\cdot\xi'-f_a(x,y)\\
&=xH_a(x)-\vfi_a(x)+\tfrac{|\xi'|^2}{2},
\qquad
x=H_a^{-1}\!\left(\xi_1+\tfrac{|\xi'|^2}{2}\right).
\end{aligned}
\ee
Here $H_a^{-1}$ denotes the inverse of the strictly increasing
function $H_a$.

We next determine $\tilde\Omega_a$ explicitly.
Put
\[
\ell_a:=\lim_{x\to\infty}H_a(x)\in\bbr\cup\{+\infty\}.
\]
Since $\tfrac{1}{a}\leq-2$, the integral
$\int_0^1(1+t)^{n-1}t^{\tfrac{1}{a}}\,dt$ diverges. Hence
\[
\lim_{x\to0^+}H_a(x)
=-\int_0^1(1+t)^{n-1}t^{\tfrac{1}{a}}\,dt
=-\infty.
\]
Together with $H_a'>0$, this shows that $H_a$ maps
$(0,\infty)$ diffeomorphically onto $(-\infty,\ell_a)$.
It follows from (\ref{4.10'}) that
\be\label{4.12}
\tilde\Omega_a
=\left\{\xi\in\bbr^n:
\xi_1+\tfrac{|\xi'|^2}{2}<\ell_a\right\}.
\ee
If $\ell_a<+\infty$, the set in (\ref{4.12}) is an open
convex sublevel set of the continuous convex function
$\xi\mapsto\xi_1+\frac{|\xi'|^2}{2}$.
If $\ell_a=+\infty$, it is all of $\bbr^n$.
Thus $\tilde\Omega_a$ is an open convex domain.
Finally, the Legendre Hessian identity gives
\[
(D_\xi^2u_b)(\xi)
=\bigl(D_{(x,y)}^2f_a(x,y)\bigr)^{-1}>0,
\qquad
\xi=\mathcal L_{f_a}(x,y).
\]
Consequently, $u_b$ is smooth and strictly convex on
$\tilde\Omega_a$.

\textbf{Proof of Theorem \ref{1.0.2}.}
Proposition \ref{4.0.1} gives the first family for
$a\in[-\tfrac12,-\tfrac1{n+1}]$.
For $b\in[-\tfrac n{n+1},-\tfrac12]$, set $a=-1-b$ and
use the dual potential $u_b$ constructed above.
By Lemma \ref{2.0.3}, it satisfies (\ref{1.1}) with exponent
$b$ and has complete Calabi metric.
Non-quadraticity is immediate from Legendre duality.
The two families cover the exponent interval
$[-\tfrac n{n+1},-\tfrac1{n+1}]$, completing the proof.
\hfill $\Box$

\begin{rmk}\label{4.0.3}
The examples above exist in every dimension $n\geq2$.
However, the rigidity proof in Section 3 relies on tools specific
to dimension two and does not extend directly to higher dimensions.
\end{rmk}

%%%%%%%%%%%%%%%%%%%%%%%
%%%%%%%%%%%%%%%%%%%%%%%%
%%%%%%%%%%%%%%%%%%%%%%%%%%%%
\vskip3mm
\paragraph*{\textbf{Funding}}
This work was supported by the National Natural Science Foundation of China (Grant No.~12501063), the NSF of Henan Province (Grant No.~ 262300421852), the Postdoctoral Fellowship Program of CPSF (Grant No.~GZC20252036), and China Postdoctoral Science Foundation (Grant No.
2025M783098).

\vskip3mm
\noindent{\bf AI assistance statement.}  The authors used AI tools only for routine drafting and computational assistance; 
the authors alone are responsible for the correctness of Theorems 1.2 and 1.3.

\vskip3mm
\paragraph*{\textbf{Conflicts of Interest}}
 \vskip2mm

The authors declare that they have no conflict of interest.

%%%%%%%%%%%%%%%%%%%%%%%%%%%%
%%%%%%%%%%%%%%%%%%%%%%%%%%%%

\end{document}